\documentclass{IEEEtran4PSCC}

\usepackage[cmex10]{amsmath}
\usepackage{amssymb}

\usepackage[hyphens]{url}
\usepackage{hyperref}
\usepackage{xurl}
\hypersetup{breaklinks=true}  % Make hyperlinks breakable

\usepackage{cite}
\usepackage{amsmath,amssymb,amsfonts,bm, bbm}
\usepackage{algorithm, algorithmic}
\usepackage{xcolor}
\usepackage[normalem]{ulem}
\usepackage{bm}
\usepackage{amsmath,amssymb,amsfonts}
\usepackage[c1]{optidef}
\usepackage{comment}
\usepackage{hyperref}
\usepackage{booktabs}
\usepackage{makecell}
\usepackage[dvipsnames]{xcolor}
\usepackage{xcolor}
\usepackage[caption=false,font=footnotesize]{subfig}
\usepackage{tabularx}
\usepackage[normalem]{ulem}
\usepackage{comment}
\usepackage{xcolor, soul}
\usepackage{tikz}
\usetikzlibrary{shapes,arrows}
\usepackage{smartdiagram}
\usepackage{multirow}

\newcommand{\black}[1]{\textcolor{black}{#1}}

\makeatletter
\let\old@ps@headings\ps@headings
\let\old@ps@IEEEtitlepagestyle\ps@IEEEtitlepagestyle
\def\psccfooter#1{%
    \def\ps@headings{%
        \old@ps@headings%
        \def\@oddfoot{\strut\hfill#1\hfill\strut}%
        \def\@evenfoot{\strut\hfill#1\hfill\strut}%
    }%
    \def\ps@IEEEtitlepagestyle{%
        \old@ps@IEEEtitlepagestyle%
        \def\@oddfoot{\strut\hfill#1\hfill\strut}%
        \def\@evenfoot{\strut\hfill#1\hfill\strut}%
    }%
    \ps@headings%
}
\makeatother

\begin{document}
\title{Constructing Deployment Scenarios for Reserve Deliverability via Adaptive Robust Optimization}

%% To specify the authors when (number of affiliations <= 2)
\author{
\IEEEauthorblockN{Guillaume Van Caelenberg$^{*}$, Akylas Stratigakos$^{\dag}$, Elina Spyrou$^{*}$}
\IEEEauthorblockA{$^{*}$Department of Electrical and Electronic Engineering, Imperial College London, U.K. \\
$^{\dag}$UCL Energy Institute, University College London, U.K.}
}

% Vectors
\def\vzero{{\mathbf{0}}}
\def\vone{{\mathbf{1}}}
\def\vmu{{\boldsymbol{\mu}}}
\def\vnu{{\boldsymbol{\nu}}}
\def\vupsilon{{\boldsymbol{\upsilon}}}
\def\vtheta{{\boldsymbol{\theta}}}
\def\vell{{\boldsymbol{\ell}}}
\def\vxi{{\boldsymbol{\xi}}}
\def\vgamma{{\boldsymbol{\gamma}}}
\def\vzeta{{\boldsymbol{\zeta}}}
\def\vpi{{\boldsymbol{\pi}}}
\def\mxi{{\boldsymbol{\Xi}}}

\def\va{{\mathbf{a}}}
\def\vb{{\mathbf{b}}}
\def\vc{{\mathbf{c}}}
\def\vd{{\mathbf{d}}}
\def\ve{{\mathbf{e}}}
\def\vf{{\mathbf{f}}}
\def\vg{{\mathbf{g}}}
\def\vh{{\mathbf{h}}}
\def\vi{{\mathbf{i}}}
\def\vj{{\mathbf{j}}}
\def\vk{{\mathbf{k}}}
\def\vl{{\mathbf{l}}}
\def\vm{{\mathbf{m}}}
\def\vn{{\mathbf{n}}}
\def\vo{{\mathbf{o}}}
\def\vp{{\mathbf{p}}}
\def\vq{{\mathbf{q}}}
\def\vr{{\mathbf{r}}}
\def\vs{{\mathbf{s}}}
\def\vt{{\mathbf{t}}}
\def\vu{{\mathbf{u}}}
\def\vv{{\mathbf{v}}}
\def\vw{{\mathbf{w}}}
\def\vx{{\mathbf{x}}}
\def\vy{{\mathbf{y}}}
\def\vz{{\mathbf{z}}}

% Matrix
\def\mA{{\mathbf{A}}}
\def\mB{{\mathbf{B}}}
\def\mC{{\mathbf{C}}}
\def\mD{{\mathbf{D}}}
\def\mE{{\mathbf{E}}}
\def\mF{{\mathbf{F}}}
\def\mG{{\mathbf{G}}}
\def\mH{{\mathbf{H}}}
\def\mI{{\mathbf{I}}}
\def\mJ{{\mathbf{J}}}
\def\mK{{\mathbf{K}}}
\def\mL{{\mathbf{L}}}
\def\mM{{\mathbf{M}}}
\def\mN{{\mathbf{N}}}
\def\mO{{\mathbf{O}}}
\def\mP{{\mathbf{P}}}
\def\mQ{{\mathbf{Q}}}
\def\mR{{\mathbf{R}}}
\def\mS{{\mathbf{S}}}
\def\mT{{\mathbf{T}}}
\def\mU{{\mathbf{U}}}
\def\mV{{\mathbf{V}}}
\def\mW{{\mathbf{W}}}
\def\mX{{\mathbf{X}}}
\def\mY{{\mathbf{Y}}}
\def\mZ{{\mathbf{Z}}}

% make the title area
\maketitle

\begin{abstract}
Network congestion often hinders the deployment of reserves needed to balance forecast errors during real-time operations.
A pertinent idea to tackle this challenge involves adding deployment scenarios of spatial distributions of forecast errors as contingencies to the day-ahead problem.
However, current approaches disregard the effect of grid characteristics and the day-ahead schedule on the induced congestion and, consequently, reserve deliverability.
In this work, we formulate a two-stage adaptive robust optimization problem to jointly consider interactions between day-ahead and real-time operations and forecast errors.
Using a column-and-constraint algorithm, we iteratively construct deployment scenarios by finding the worst-case forecast error for reserve deliverability.
Simulations on the RTS-GMLC system show that adding these scenarios to the day-ahead problem significantly reduces the frequency of \black{congestion-driven reserve undeliverability.}
Notably, \black{the choice and number of scenarios dynamically adapts to the day-ahead schedule}.
\end{abstract}

\begin{IEEEkeywords}
Electricity markets, \black{operational uncertainty}, probabilistic forecasting, adaptive robust optimization, \black{reserve deliverability}
\end{IEEEkeywords}

\thanksto{\noindent This work was funded by the Leverhulme International Professorship (LIP-2020-002) and the Engineering and Physical Sciences Research Council under the grant EP/Y025946/1 (Electric Power Innovation for a Carbon-free Society (EPICS)).}

\section{Introduction} 
\label{introduction}
\paragraph{Background and motivation}
To manage the growing imbalances between the day-ahead (DA) and real-time (RT) scheduling problems, partly caused by increasing variable renewable energy (VRE) forecast errors, 
system operators (SOs) reserve spare generation capacity in DA, 
which can be deployed in RT \cite{ferc2021}.
To ensure that the frequency of reserve shortages remains below a prescribed reliability target, 
SOs are increasingly using probabilistic forecasts to determine the amount of required reserves dynamically 
\cite{devos2019dynamic, haupt2019_use_prob_forecast}.
Typically, reserve requirements are estimated for the system as a whole or a limited number of large zones \cite{epri2019_ancillary_services}. 
However, network congestion often renders reserves undeliverable in RT. 
For instance, a significant portion of ramping reserve was undeliverable in 2018 in California \cite{caiso2020frp_training}.

\paragraph{Related work}
Several approaches of varying complexity could address the problem of undeliverable reserves.
At one end of the spectrum, computationally intensive formulations explicitly model multiple scenarios for the spatial distribution of forecast errors to obtain DA schedules that ensure the desired RT reliability \cite{roald2023power}. 
For instance, two-stage stochastic optimization problems endogenously determine reserve capacity at various locations by minimizing probability-weighted costs over multiple scenarios with distinct spatial distributions of forecast errors. 
However, such methods may be impractical due to computational issues and \black{lack of properties, desirable in markets \cite{kazempour2018astochastic}}.
As an alternative, robust optimization problems \cite{bertsimas2020_rob_and_adapt_opt} ensure deliverability for all error realizations within an \textit{uncertainty set} \black{and closely resemble existing scheduling processes \cite{street2025robustness}}.
\black{However, robust problems are sensitive to the design of the uncertainty set and often lead to overly conservative solutions.} 
\black{This conservativeness is to some extent reduced when robust problems are adaptive,  considering multiple stages and allowing recourse actions at a later stage (e.g., RT redispatch)  \cite{bertsimas2012adaptive}.}
Solution approaches for adaptive robust problems usually assume a reserve activation policy \cite{bertsimas2024marginal} or iteratively find the worst-case realizations for a series of DA schedules  \cite{donti2021_adversarially}. 
\black{Data-driven methods can be used to determine the parameters of reserve activation policies \cite{lyon2015_locational_disqualification, singhal2019data_driven_reserve_response}.}

At the other end, less computationally intensive methods rely on heuristics and approximations for generating DA schedules. 
\black{While they aim for higher RT reliability than practices with network-unaware reserve procurement, they usually lack theoretical guarantees.}
For instance, reserving capacity at each node to meet nodal reserve needs \black{is an inner approximation of the adaptive robust problem's feasible region that does not require any network capacity for reserves}.
In practice, imposing reserve \black{requirements} at every node is often costly and sometimes impossible due to the lack of local reserve suppliers.
\black{When reserve requirements are nodal, there are provisions for limited contribution of out-of-node resources. 
The limits on the contributions are determined through conservative approaches \cite{madani2025_inscribed_box} or heuristics \cite{liang2022_weather_driven, vandenbergh2018cross_border}.}

\black{Building upon the current practice of approximations with zonal reserve requirements, 
several articles propose a shift from zones with jurisdictional boundaries \cite{meeus2020_evolution_electricity_markets}} to
 zones identified based on active and reactive power flow analysis \cite{kumar2004_zonal_congestion_management} and grid characteristics \cite{wang2015, xu2016_zonal_ordc}.
\black{
To account for dynamic congestion patterns,  \cite{viafora022dynamic_reserve} presented a method for dynamically partitioning the grid into zones using probabilistic forecasts.}

More recently, SOs started incorporating forecast error scenarios as contingencies in their DA schedule to support reserve deliverability in RT \cite{caiso2020frp, chen2023_addressing_uncertainties};
hereafter, we refer to them as \textit{deployment scenarios}.
By explicitly modeling RT redispatch actions under deployment scenarios, this approach advances the traditional deterministic modeling practice\cite{hobbs2022can} and constitutes a second-best alternative to the optimal, but impractical, multi-stage stochastic and robust formulations.
However, generating these deployment scenarios from probabilistic forecasts is an open question. 
Industry best practices typically account for a few \textit{extreme scenarios} with heuristically selected spatial distributions of forecast errors \cite{caiso2025_dame_business_requirements}, 
potentially overlooking other possible distributions that could lead to reserve undeliverability due to network topology and transmission constraints.
Importantly, these deployment scenarios are chosen without considering how reserve deliverability depends on the DA schedule and network congestion.

\paragraph{Aim and contribution}
In this work, we develop a method for constructing deployment scenarios from a set of probabilistic forecasts, which are then included as contingencies in the DA scheduling problem.
We formulate a two-stage adaptive robust optimization (ARO) problem that jointly considers the DA schedule, 
worst-case forecast error realizations, and the full-recourse RT schedule.
When the uncertainty set is polyhedral, the ARO problem is equivalent to the standard DA scheduling problem with deployment scenarios as additional contingencies.
Using the column-and-constraint (CCG) \cite{zeng2013solving_two_stage} algorithm, we construct deployment scenarios corresponding to a series of iteratively updated DA schedules.
We apply this method to the RTS-GMLC 2019 \cite{rts_gmlc} system and provide original insights into the differences between \textit{CCG-constructed} deployment scenarios and \textit{extreme} deployment scenarios constructed following industry practice.
We find that CCG-constructed deployment scenarios often include errors in both directions (i.e. actuals may be greater or lower than the DA point forecast) at different grid nodes,
which contrasts {extreme scenarios} that include errors in a single direction. 
To account for the aggregation benefits of forecast errors, extreme scenarios reduce the reserve requirement at each node based on a uniform across-nodes policy.
In contrast, the CCG-constructed deployment scenarios allow the DA problem to endogenously calculate what level of aggregation benefit (if any) can be leveraged, considering transmission constraints.
Overall, our results show that incorporating deployment scenarios in the DA schedule improves reserve deliverability compared to network-unaware reserve procurement and that CCG-constructed deployment scenarios closely approximate prescribed reliability targets, 
while outperforming extreme scenario approaches.

The rest of the paper is structured as follows. 
Section~\ref{section_prelim} presents the preliminaries on the reserve deliverability problem. 
Section~\ref{section_method} describes the methodology proposed to construct deployment scenarios. 
Section~\ref{section_results} presents the experimental design and results. 
Finally, Section~\ref{section_conclusions} summarizes conclusions and provides directions for further research.

\section{Preliminaries}\label{section_prelim}
In this section, we describe the operating framework (Subsection~\ref{method_preliminaries}), 
and formulate the DA scheduling problem with deployment scenarios for the deliverability of reserves as a two-stage ARO problem (Subsection~\ref{aro_subsection}).

\paragraph*{Notation} 
We use bold lowercase (uppercase) font for vectors (matrices) and calligraphic font for sets. 
Let $|\cdot|$ be the set cardinality and $\vone$ be a vector of ones with appropriate size.
Forecasts are denoted by $\hat{(\cdot)}$;
the value of decision variables at optimality is denoted by $\dot{(\cdot)}$.
\black{For a positive integer $Q$, we define $[Q] = \{1,\dots, Q\}$.}

\subsection{Operating framework}\label{method_preliminaries}

Consider a power system where $\mathcal{N}$ is the set of nodes, $\mathcal{L}$ is the set of lines, and $\mathcal{G}$ is the set of generators.
For the $i$th period, the point forecast for the net demand (load minus VRE production) is $\hat\vd\in\mathbb{R}^{|\mathcal{N}|}$, and $\vxi\in\mathbb{R}^{|\mathcal{N}|}$ is a random variable that describes the forecast error, 
defined as actual minus DA point forecast. 
The SO has access to a set \black{$\mathcal{E}$ of $K$ scenarios with forecast errors that account for spatial correlations, where} $\mathcal{E} = \{\hat\vxi_{k}\}_{k\in[K]}$.
In practice, this information can be provided by methods that generate multivariate scenarios from marginal probabilistic forecasts using a Copula function --- see, e.g.,\cite{carmona2024_joint_granular_model}. 
We note the aggregate net demand error at the system level with $\black{\xi^{\textrm{agg}}_{}}$ \black{and compute its scenario-specific value, $\hat\xi^{\textrm{agg}}_{k}$}, as $ \vone^{\top}\ \hat\vxi_{k}$.
For brevity, the period index $i$ is omitted.

Consider a prescribed reliability level $\alpha \in (0,1)$, where $\alpha$ is typically high, e.g., $0.90$, $0.95$. 
The system-wide reserve requirements for the $i$th period are estimated by
\begin{align}\label{sw_reserve_req}
\black{\hat \rho^{+} =  q_{\frac{1+\alpha}{2}}\big(\{\hat\xi^{\textrm{agg}}_k\}_{k\in[K]}\big),
\,\hat \rho^{-} =  q_{\frac{1-\alpha}{2}}\big(\{\hat\xi^{\textrm{agg}}_k\}_{k\in[K]}\big)}
\end{align}
where \black{$q_{u}(\cdot)$ is a function that takes as input $K$ scenarios and returns the $u$th empirical quantile}.
In words, the likelihood of $\xi^{\textrm{agg}}$ being in the interval $[\hat\rho^{-}, \hat\rho^{+}]$ is $100\cdot\alpha\%$.
Typically, it holds that $\hat\rho^{-}\leq0$ and $\hat\rho^{+}\geq0$, 
\black{which implies that downward and upward reserves are procured to balance negative and positive errors of net demand, respectively.}

The SO solves two sequential optimization problems. 
In DA, the SO schedules energy and reserves. 
In RT, the SO adjusts the schedules to ensure cost-effective and reliable operations. 
As we are concerned with reserve undeliverability, the second stage objective only considers reliability. 
\footnote{\black{
To highlight the article's novel contributions, we solve a single-interval optimization problem.
The core methodology can be extended for scheduling applications with commitment decisions and inter-temporal constraints.}}

\paragraph{DA problem}
The DA (or first-stage) scheduling problem is given by 
\begin{subequations}\label{dsw} 
\begin{align} 
\underset{\substack{\vp, \vr^+, \vr^-}}{\min} & \,\, (\vc^{\textrm{e}})^{\top}\vp  + (\vc^{+})^{\top}\vr^{+} + (\vc^{-})^{\top}\vr^{-}, \label{dsw_obj} & \\
\text{s.t.} & \, \, \vone^{\top}\vp = \vone^{\top}\hat\vd, &  \label{dsw_c1} \\
& \, \, \vone^{\top}\vr^+ \geq \hat{\rho}^{+}, \vone^\top \vr^- \geq -\hat{\rho}^{-}, & \label{dsw_c2} \\
& \, \, \black{\mM(\mA\vp - \hat\vd) \leq \vf^{\max}}, & \label{dsw_c3}\\
& \, \, \black{-\mM(\mA\vp - \hat\vd) \leq \vf^{\max}}, & \label{dsw_c4}\\
& \, \, \vp + \vr^+ \leq \vp^{\max}, \label{dsw_c5}\\
& \,\, \vp - \vr^- \geq \vp^{\min},  \label{dsw_c6}\\ 
& \, \, \vp, \vr^+, \vr^- \geq \vzero.  & \label{dsw_c7}
\end{align}
\end{subequations}
Problem \eqref{dsw} computes the least-cost
energy ($\vp\in \mathbb{R}^{|\mathcal{G}|}$) and reserve schedule  ($\vr^+, \vr^-\in \mathbb{R}^{|\mathcal{G}|}$) \eqref{dsw_obj} 
to satisfy the net demand forecast \eqref{dsw_c1} and dynamic reserve requirements \eqref{dsw_c2}, 
subject to transmission feasibility constraints \eqref{dsw_c3}-\eqref{dsw_c4} and constraints for minimum and maximum generating power \eqref{dsw_c5}-\eqref{dsw_c7}.
Here, $(\vc^{\textrm{e}}, \vc^{+}, \vc^{-})$ represent linear costs,
$(\vp^{\textrm{max}}, \vp^{\min}, \vf^{\textrm{max}})$ represent technical limits for energy and line flows, 
\black{$\mM \in {\mathbb{R}}^{|\mathcal{L}| \times |\mathcal{N}|}$ is the Power Transfer Distribution Factors (PTDF) matrix, and $\mA \in {\mathbb{R}}^{|\mathcal{N}| \times |\mathcal{G}|}$ is the node incidence matrix that maps generators to nodes.}

\paragraph{RT problem}
Given fixed first-stage decisions $\dot{\vx}=(\dot\vp, \dot\vr^{+}, \dot\vr^{-})$ and a realized net demand forecast error $\vxi_0$, the RT (or second-stage) redispatch problem that computes recourse actions $\vp^{\textrm{rec}}\in{\mathbb{R}}^{|\mathcal{G}|}$ is given by
\begin{subequations}\label{rt-prob} 
\begin{align} 
\underset{\substack{\vp^{\textrm{rec}}, 
\vg^{+}, \vg^{-}, \\\vell^{+}, \vell^{-}
}}{\min} & \,\, (\vc^{\textrm{viol}})^\top (\vg^{+}+\vg^{-}+\vell^{+}+\vell^{-}), \label{rt_obj} & \\
\text{s.t.} & \, \,  \vone{^\top}\vp^{\textrm{rec}} = \vone{^\top} \vxi_0, & \label{rt_c1}
\\
& \, \, \black{\mM(\mA} \vp^{\textrm{rec}} - \vxi_0) \leq \black{\vf^{\max} - \dot\vf^{e}+\vell^+}, & \label{rt_c2}\\
& \, \, -\black{\mM(\mA}\vp^{\textrm{rec}} - \vxi_0) \leq\black{\vf^{\max} + \dot\vf^{e}+\vell^-}, & \label{rt_c3}\\
& \, \, \vp^{\textrm{rec}} \leq \dot\vr^{+} +\vg^+, & \label{rt_c4}
\\ & \, \, -\vp^{\textrm{rec}} 
\leq \dot\vr^{-} +\vg^- , & \label{rt_c5}
\\ & \, \, \vg^+, \vg^-, \vell^+, \vell^- \geq \vzero, & \label{rt_c6}
\end{align}
\end{subequations}
\black{where $\dot\vf^{e}=\mM(\mA\dot\vp - \hat\vd)$ is shorthand for the DA scheduled energy flow. 
Note that the available transmission capacity for recourse actions is $\vf^{\max} - \dot\vf^{e}$ and $\vf^{\max} + \dot\vf^{e}$.} 
Constraint \eqref{rt_c1} ensures the RT balance of supply and demand.
Constraints \eqref{rt_c2}-\eqref{rt_c3} and \eqref{rt_c4}-\eqref{rt_c5} ensure that recourse actions are within the available transmission capacity and generation reserved in DA, respectively. 
\black{If the reserves are insufficient, e.g., RT schedules exceed capacities reserved in DA or cannot be delivered due to grid congestion, }
then the slack variables $\vg^-, \vg^+ \in \mathbb{R}^{|\mathcal{G}|}$  
and $\vell^-, \vell^+ \in \mathbb{R}^{|\mathcal{L}|}$ become positive \cite{wang2020market}. 
The objective function \eqref{rt_obj} imposes a large penalty $\vc^{\textrm{viol}}$ on non-zero slack values.

\subsection{Adaptive Robust Optimization Formulation}\label{aro_subsection}

We now formulate a two-stage ARO problem that enables reserve deliverability 
and show that it can be equivalent to \eqref{dsw} augmented with deployment scenarios.
To streamline notation, we rewrite problem \eqref{dsw} in a compact form as
\label{dsw-std} 
\begin{align*} 
\underset{\substack{\vx}}{\min} & \,\, \vc^{\top}\vx, \,\, \text{s.t.} \, \, \black{\mB \vx \leq \vb},  \label{} 
\end{align*}
where $\vx =(\vp, \vr^+, \vr^-)$, $\vc=(\vc^e, \vc^+, \vc^-)$, \black{and $(\mB, \vb)$ parameterize the set of linear inequalities\footnote{\black{Each equality constraint is replaced with two opposite inequalities.}} that represent the feasible set of \eqref{dsw}.} 
Next, we rewrite problem \eqref{rt-prob} as
\label{rt-std} 
\begin{align*} 
\underset{\substack{\black{\vp^{\textrm{rec}}},\vs\geq \vzero}}{\min} & \,\, (\vc^{\textrm{viol}})^{\top}\vs, & \text{s.t.} & \, \, \black{\mH\dot \vx +\mD \vp^{\textrm{rec}} - \vs \leq \mE \vxi_0 + \vh}, &  \label{} 
\end{align*}
where $\dot\vx$ are the fixed first-stage decisions, $\vs=(\vzero, \vg^+, \vg^-, \vell^+, \vell^-)$,
 and $\black{(\mH, \mD, \mE, \vh)}$ are constructed appropriately to represent the feasible set of \eqref{rt-prob}.\footnote{
 \black{In addition to slack variables, $\vs$ includes  $\vzero$ for the two opposite inequalities, replacing  \eqref{rt_c1}}.}

The two-stage ARO problem with full recourse, considering only feasibility penalties, is given by
\begin{subequations}\label{2ro-prob} 
\begin{align} 
\underset{\black{\substack{\vx}}}{\min} & \,\, \vc^{\top}\vx + \underset{\substack{\vxi\in\mathcal{U}}}{\max} 
\,\, \underset{\substack{\vp^{\textrm{rec}}(\vxi), \vs(\vxi)}}{\min} \,\, (\vc^{\textrm{viol}})^{\top}\vs(\vxi), \label{} & \\
\text{s.t.} & \, \, \black{ \mB\vx \leq \vb}, &  \label{} \\ 
& \, \, \black{ \mH\vx + \mD\vp^\textrm{rec}(\vxi) {- \vs(\vxi)} \leq \mE\vxi + \vh}, & \forall \vxi \in {\mathcal{U}}, &  \label{} \\
& \,\, \vs(\vxi) \geq \vzero, & \forall \vxi \in {\mathcal{U}},
\end{align}
\end{subequations}
where second-stage decisions $\vp^\textrm{rec}(\vxi), \vs(\vxi)$ are a function of $\vxi$
and $\mathcal{U}$ is an uncertainty set that covers potential realizations of net demand error $\vxi$ \black{the system has to be reliable against. For example, $\mathcal{U}$ could coincide with $\mathcal{E}$, or be a subset or an approximation of it}.

The choice of $\mathcal{U}$ is critical for both the performance guarantees and the computational cost associated with solving (or approximating) \eqref{2ro-prob}.
When $\mathcal{U}$ is a polyhedral set, 
then the worst-case cost occurs at one of its extreme points (vertices).
Let $\mathcal{V}$ be a discrete set that contains all the vertices of $\mathcal{U}$.
Then, \eqref{2ro-prob} is equivalent to 
\begin{subequations}\label{2ro-reform} 
\begin{align} 
\underset{\substack{\vx, \{\vp^{\textrm{rec}}_\vxi, \vs_{\vxi}\}, \eta } }{\min} 
& \,\, \vc^\top\vx + \eta, \label{2ro-reform-obj} &
\\
\text{s.t.} & \, \, \black{\mB\vx \leq \vb,} &\\
& \, \, \eta \geq (\vc^{\textrm{viol}})^\top\vs_\vxi, & \forall \vxi \in {\mathcal{V}}, &\\
& \, \, (\vp^{\textrm{rec}}_\vxi, \vs_{\vxi}) \in \Omega(\vx, \vxi), & \forall \vxi \in {\mathcal{V}}, &\label{c3_2ro-reform}
\end{align}
\end{subequations}
where $\eta$ represents the worst-case violation cost, $(\vp^{\textrm{rec}}_\vxi, \vs_{\vxi})$ are the wait-and-see second-stage decisions for each vertex in $\mathcal{V}$, 
and 
$\Omega(\vx, \vxi) = \{ (\vp^{\textrm{rec}}, \vs): \black{\mH\vx + \mD\vp^{\textrm{rec}}(\vxi) - \vs(\vxi) \leq \mE\vxi + \vh,
\,\vs(\vxi)\geq \vzero }\}$ is the feasible set of the second-stage problem, given $\vx, \vxi$.
The objective function \eqref{2ro-reform-obj} consists of the DA scheduling cost and the worst-case violation cost. 
\black{Note that \eqref{c3_2ro-reform} includes the RT constraints for particular realizations of $\vxi$, which could be interpreted as reserve deployment scenarios.
Hence, the ARO formulation formally justifies the industry intuition of adding deployment scenarios to the DA problem.}
When \eqref{c3_2ro-reform} is \black{satisfied without slack activation for all vertices, $\eta$ equals zero.}
\black{If any slack is non-zero, 
then $\eta$ is positive, and supply-demand balance is not guaranteed within $\mathcal{U}$.}
\black{In those cases, it is worth examining whether the system lacks reserves at the system level or at specific nodes and how accurate model inputs, such as $\mathcal{U}$  and $\vc^{\textrm{viol}}$, are.}

The ARO reformulation shows that the success of the deployment scenarios hinges on the careful construction of these scenarios.
The number of vertices in $\mathcal{U}$ may be exponential to the dimension of $\vxi$, 
making vertex enumeration impractical for solving \eqref{2ro-reform} to optimality. 

\section{Constructing Deployment Scenarios}\label{section_method}

In this section, we \black{describe two methods for constructing deployment scenarios, assuming a polyhedral uncertainty set presented in Subsection~\ref{set_subsection}.} 
First, we describe a method that only uses DA forecasts of net demand as inputs and resembles industry practices (Subsection~\ref{for_depl_subsection}). 
Second, we develop a method that additionally leverages grid characteristics and DA decisions, which comprise our key contribution (Subsection~\ref{ccg_subsection}).

\subsection{Uncertainty Set}\label{set_subsection}

In this work, we assume that SOs aim to ensure reserve deliverability for RT scenarios with $\xi^{\textrm{agg}}$ lying within an interval bounded by the \black{system-wide} reserve requirements $\hat\rho^+, \hat\rho^-$, and nodal net demand errors lying within a ``box'', given by
\begin{align}
&
\mathcal{U}^{\textrm{agg}}_{\alpha} = \{ \vxi \,:\,  \hat\rho^{-} \leq  \vone^\top \vxi \leq \hat\rho^{+}\},
\\        
& \mathcal{U}^{\textrm{box}} = \{ \vxi \,:\, \min_{k\in[K]}{\{\hat\vxi_k}\} \leq \vxi \leq
        \max_{k\in[K]}{\{\hat\vxi_k}
        \}\},\label{box_set}
\end{align}
where \black{the box is bounded by the lowest and highest forecast error scenario (element-wise operation).}

The uncertainty set for problem \eqref{2ro-prob} is given by the intersection $\mathcal{U}=\mathcal{U}^{\textrm{agg}}_{\alpha}\cap\mathcal{U}^{\textrm{box}}_{}$, 
which implies that $\mathcal{U}\subseteq\mathcal{U}^{\textrm{agg}}_{\alpha}$ and $\mathcal{U}\subseteq\mathcal{U}^{\textrm{box}}_{}$.
\black{In practice, SOs can also use other polyhedral uncertainty sets that may better reflect the distribution of forecast errors. The quality of uncertainty sets can be assessed in terms of coverage, density, or interpretability.}

\subsection{Extreme Deployment Scenarios from Forecasts}\label{for_depl_subsection}

This approach resembles industry practice \cite{caiso2025_dame_business_requirements} and constructs a set of two extreme deployment scenarios $\mathcal{S}^{\textrm{depl}}$: one for positive \black{and one for negative forecast errors.}
This approach assumes that the most challenging case for reserve deliverability is when all errors are in the same direction (positive or negative), which implies a high correlation among them.
Here, we follow the same logic and distribute $\hat{\rho}^{+}, \hat{\rho}^{-}$ to the nodes such that all errors are in the same direction. First, for the $j$th node, 
we calculate allocation factor $e^+_j$, which accounts for the $\frac{1+\alpha}{2}$-level marginal quantile forecasts as follows:
\begin{align}\label{allocation_factors}
\black{    e^+_j = \frac{ q_{\frac{1+\alpha}{2}}\big(\{\hat\xi_{j,k}\}_{k\in[K]}\big)}
    {\sum_{n\in\mathcal{N}} q_{\frac{1+\alpha}{2}}\big(\{\hat\xi_{n,k}\}_{k\in[K]}\big)}.}
\end{align}
\black{In words, the allocation factors are proportional to the quantiles of nodal forecast errors.}
We construct the upward deployment scenario that includes in each node $n$ an upward error equal to $\xi_n^+=\hat\rho^+\cdot e^+_n$.
\black{An additional Euclidean projection step onto $\mathcal{U}$ is applied.}
We follow the same approach for the downward direction.

\subsection{Constructing Deployment Scenarios via CCG}\label{ccg_subsection}

In this section, we construct deployment scenarios considering their dependency on \emph{(i)} problem parameters (e.g., grid topology) and \emph{(ii)} first-stage decisions.
\black{We apply a CCG algorithm \cite{zeng2013solving_two_stage} to iteratively add vertices, such that $\mathcal{S}^{\textrm{depl}}\subseteq\mathcal{V}$}.

To streamline notation, let
\begin{align}\label{rt_valuefunc}
   \black{{Q}(\dot\vx) = \underset{\substack{\vxi\in\mathcal{U}}}{\max} 
\,\, \underset{\substack{\vp^{\textrm{rec}}(\vxi), \vs(\vxi)\in \Omega}(\dot\vx, \vxi)}{\min} \,\, (\vc^{\textrm{viol}})^\top\vs(\vxi)}
\end{align}
be the worst-case objective value of the second-stage problem given fixed first-stage decisions $\dot\vx$.
The CCG algorithm is summarized in Algorithm~\ref{CCG_algo}. 
{First, we consider an empty set of deployment scenarios,} \black{$\mathcal{S}^{\textrm{depl}}$}, 
solve \eqref{2ro-reform},
and estimate \black{a lower bound $\texttt{LB}$ on the RT violation cost} (lines~\ref{ccg_step1}, \ref{ccg_step3}).
Next, we fix the first-stage decisions and approximate the worst-case cost of the second stage $Q(\dot\vx)$ (line~\ref{ccg_step4}), 
which returns an approximately worst-case scenario $\tilde\vxi$ and an upper bound \texttt{UB} \black{on the RT violation cost }(line~\ref{ccg_step5}).
If the worst-case scenario $\tilde\vxi$ \black{has a non-zero violation cost}, 
it is added to the set of deployment scenarios (lines~\ref{ccg_step6}-\ref{ccg_step9}).
\black{The algorithm terminates when $\texttt{UB}=\texttt{LB}$ or when the maximum number of iterations $M$ is reached. $M$ can be chosen based on the number of deployment scenarios the SO can add to the DA schedule.}

{If an oracle that optimally solves ${Q}(\dot\vx_j)$ for a fixed $\dot\vx_j$ is available (step~\ref{ccg_step4}), 
then \eqref{2ro-reform} is solved to optimality when Algorithm~\ref{CCG_algo} converges.}
In practice, however, computing ${Q}(\dot\vx_j)$ is challenging because problem~\eqref{rt_valuefunc} is a nonconvex $\max\textrm{-}\min$ problem.
Solution approaches reformulate the problem as a mixed-integer problem or rely on heuristics\cite{zeng2013solving_two_stage}.
Here, we use the \textit{alternating direction method} (ADM) heuristic to approximate ${Q}(\dot\vx_j)$ \eqref{rt_valuefunc}, which has been shown to perform well in similar problems \cite{lorca2014adaptive}.

First, strong duality is applied to the inner $\min$ problem for a fixed $\dot\vx$, reformulating the $\max\textrm{-}\min$ problem into a $\max$ problem given by
\begin{subequations}\label{adm-inner}
\begin{align}
\underset{\substack{\vpi, \vxi \in \mathcal{U}} }{\max} 
& \,\, \black{\vpi^\top(\mH\dot\vx-\mE\vxi - \vh)}, \label{} & \\
\text{s.t.} & \, \, \black{-\mD^\top\vpi = \vzero},
\\ & \, \, \vzero \leq \vpi \leq \vc^{\textrm{viol}}, & \label{adm-inner-c2}
\end{align}    
\end{subequations}
where $\vpi$ denotes dual variables of the constraints in $\Omega(\dot\vx, \vxi)$ \footnote{We slightly abuse notation as \eqref{adm-inner-c2} should be omitted for the dual corresponding to \eqref{rt_c1} as it concerns an equality constraint with no feasibility slack.}.
Second, Problem \eqref{adm-inner}, which has a bilinear objective, is solved via Algorithm~\ref{ADM_algo}.
The algorithm relies on iteratively optimizing a relaxed linear program, where part of the decision variables are treated as constants (either $\vpi$ or $\vxi$).
Given fixed first-stage decisions $\dot\vx$ and an initial guess $\vxi^{\textrm{init}}$, we maximize \eqref{adm-inner} over $\vpi$ (assuming $\vxi=\vxi^{\textrm{init}}$), 
which provides a lower bound $\texttt{LB}^Q$ (line~\ref{adm_step3}).
Then, we maximize \eqref{adm-inner} over $\vxi\in \mathcal{U}$, while $\vpi$ is fixed at the solution, 
which provides an upper bound $\texttt{UB}^Q$ (line~\ref{adm_step4}).
\black{The local bounds are iteratively updated until convergence ($\texttt{UB}^Q = \texttt{LB}^Q$), 
which is then guaranteed to be a local optimum that satisfies the Karush–Kuhn–Tucker conditions \cite{lorca2014adaptive},
or the algorithm terminates when the maximum number of iterations $L$ is reached.}
\black{We set $\tilde{Q}(\dot\vx_j)$ as the average between $\texttt{UB}^Q, \texttt{LB}^Q$ and return the worst-case scenario $\vxi^{\textrm{wc}}=\vxi_l$ (line~\ref{adm_step7}).}

The initial guess $\vxi^{\textrm{init}}$ can be critical to the convergence of Algorithm~\ref{ADM_algo}.
\black{Here, we develop a heuristic that chooses as initial guesses forecast error scenarios that could aggravate congestion.
First, we select a set of lines that are congested in at least one period in DA or RT when the solution of \eqref{dsw} is followed.
Then, for each line $l$, we find the forecast error $\vxi$ that results in the highest power flow increase in the same direction as the DA flow $\dot f_l^{e}$, obtained from \eqref{dsw}, as follows}
\begin{align}\label{adm_init}
\underset{\substack{\vxi\in\mathcal{U}^{\textrm{box}}}}{\max}
    \,|\dot f^{e}_l - \mM_l\vxi|.
\end{align}
\black{\eqref{adm_init} admits a closed-form solution (note that $\dot f^{e}_l$ is fixed), which is a vertex of $\mathcal{U}^{\textrm{box}}$. The vertex includes, for each node $j$, the maximum nodal forecast error, when term $(\textrm{sign}(\dot f^{e}_l)\cdot M_{l,j})$ is negative; and the minimum nodal forecast error when the same term is positive.  This vertex is then projected onto $\mathcal{U}$ to obtain an initial guess.
In addition, we use the two extreme deployment scenarios from Subsection~\ref{for_depl_subsection} as initializations $\vxi^{\textrm{init}}$. That way, the use of the CCG-constructed scenarios will result in a level of reliability at least as good as the one achieved by using the extreme scenarios.
Out of all initializations, we pick the scenario $\vxi^{\textrm{wc}}$ that leads to the highest $\tilde{Q}(\dot\vx_j)$.}

\black{The computational overhead of the CCG algorithm is anticipated to be similar to that of other iterative solutions currently implemented by system operators \cite{chen2023_sc_uc_modeling}.
Note that Algorithm 2 can be trivially parallelized for multiple initializations.}

\begin{algorithm}[tb!]
\caption{ \texttt{Column-and-constraint Generation}}\label{CCG_algo}
  \textbf{Input:} Problem \eqref{2ro-reform}, maximum number of scenarios $M$.  
  \\
  \textbf{Output:} $\mathcal{S}^{\textrm{depl}}$
  \begin{algorithmic}[1]
  \STATE \label{} 
  Initialize $\mathcal{S}_0=\emptyset$, $\texttt{UB}=\infty$, $\texttt{LB}=-\infty$, $j=0$.\label{ccg_step1}
\WHILE{$\texttt{UB}-\texttt{LB}\geq 0$ and $j \leq M$}\label{ccg_step2}
\STATE Solve \eqref{2ro-reform} with $\mathcal{V}=\mathcal{S}_j$ and set $\dot\vx_j=\vx^*$, $\texttt{LB}=\eta^*$.\label{ccg_step3}
\STATE Approximate $\mathcal{Q}(\dot\vx_j)$ using $\mathcal{\tilde Q}(\dot\vx_j)$.\label{ccg_step4}
\STATE Set $\vxi_j^{\textrm{wc}}=\tilde\vxi$, 
{\texttt{UB} = $\mathcal{\tilde Q}(\dot\vx_j)$}.\label{ccg_step5}
\IF{$\mathcal{\tilde{Q}}(\dot\vx_j)>0$} \label{ccg_step6}
\STATE Update $\mathcal{S}_{j+1}=\mathcal{S}{_j} \cup \{\vxi_j^{\textrm{wc}}\}$, $j\leftarrow j+1$.
\ENDIF\label{ccg_step_8}
\ENDWHILE \label{ccg_step9}
  \STATE Return $\mathcal{S}^{\textrm{depl}} = \mathcal{S}_j$.\label{ccg_step10}
  \end{algorithmic}
\end{algorithm}

\begin{algorithm}[tb!]
\caption{ \texttt{Alternating Direction Method}}\label{ADM_algo}
  \textbf{Input:} First-stage decisions $\dot\vx$, maximum number of iterations $L$, initialization $\vxi^{\textrm{init}}$.
  \\
  \textbf{Output:} $\mathcal{\tilde Q}(\dot\vx)$, $\tilde\vxi$
  \begin{algorithmic}[1]
  \STATE \label{adm_step1} Initialize by $\texttt{UB}^Q=\infty$, $\texttt{LB}^Q=-\infty$, $l=0$, $\vxi_l=\vxi^{\textrm{init}}$.
  \WHILE{$\texttt{UB}^Q-\texttt{LB}^Q\geq \epsilon$ and $l \leq L$}\label{adm_step2}
  \STATE Solve \black{$\underset{\substack{ \vzero\leq \vpi \leq \vc^{\textrm{viol}} }}{\max} \, \vpi^\top(\mH\dot\vx - \mE\vxi - \vh), \, \text{s.t.} -\mD^\top \vpi = \vzero$},\\
  set $\vpi_l=\vpi^*$, update $\texttt{LB}^Q$. \label{adm_step3}
  \STATE Solve $\underset{\substack{\vxi \in \mathcal{U}}}{\max} \, \black{\vpi_l^\top(\mH\dot\vx-\mE\vxi-\vh),}$ update $\texttt{UB}^Q$.\label{adm_step4}
  \STATE $l\rightarrow l+1$\label{adm_step5}
  \ENDWHILE\label{adm_step6}
  \STATE Return $\mathcal{\tilde Q}(\dot\vx)=\frac{\texttt{UB}^Q+\texttt{LB}^Q}{2}$, $\tilde\vxi\leftarrow\vxi_l$.\label{adm_step7}
  \end{algorithmic}
\end{algorithm}

\section{Numerical Experiments and Results}
\label{section_results}

In this section, we summarize  the experimental design (Subsection \ref{benchmarks_subsection}), discuss an illustrative example (Subsection \ref{illustrative_ex_subsection}) and analyze results for annual simulations of the RTS-GMLC 2019 system (Subsection \ref{large_system_subsection}).

\subsection{Experimental Design} \label{benchmarks_subsection}

For a given reliability level $\alpha$, we contrast three approaches: \emph{(i)} \texttt{DSW}, where we solve \eqref{dsw} with system-wide reserve requirements estimated from \eqref{sw_reserve_req};
\emph{(ii)} $\texttt{EXT}$, where we solve \eqref{2ro-reform} with two {extreme} deployment scenarios constructed with method shown in Subsection~\ref{for_depl_subsection}; 
and \emph{(iii)} $\texttt{CCG}_{}$, 
where we solve \eqref{2ro-reform} with deployment scenarios constructed by applying  CCG, as shown in ~\ref{ccg_subsection}, 
for $\mathcal{U}=\mathcal{U}_\alpha^{\textrm{agg}}\cap\mathcal{U}^{\textrm{box}}_{}$.

\subsection{Illustrative Example} \label{illustrative_ex_subsection}

\begin{figure}[tb!]
    \centering
    \resizebox{0.75\columnwidth}{!}{
    \input{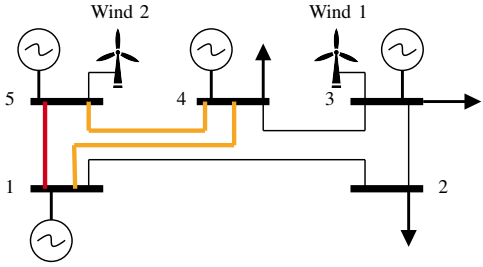}}
    \caption{Schematic of the 5-bus system.
    Red color indicates line congestion \black{in DA \texttt{DSW} schedule}.
    Orange color indicates RT line congestion given the \texttt{DSW} decisions, \black{causing deliverability problems.}}
\label{5bus_schematic}
\end{figure}

We first discuss results for a modified IEEE 5-bus system \cite{mieth2023prescribed}. 
{The results can be easily verified by the reader and provide intuition for appreciating the results for the larger system in the next section.}
The system has two wind power plants, shown in Fig.~\ref{5bus_schematic}. 
Wind plant $1$ is located in node $3$ alongside the most expensive generator, 
whereas wind plant $2$ is located in node $5$ alongside the cheapest generator.
We consider wind forecast errors following a zero-mean multivariate normal distribution, where 
\begin{equation*}
\bm{\Sigma} = \begin{bmatrix}
            0.141 & 0.001\\
            0.001 & 0.141
        \end{bmatrix} {\rm (pu)}^2,
\end{equation*}
is the covariance matrix.
We sample $K=1\,000$ DA scenario forecasts for $\hat\vxi_k$,
calculate the system-wide requirements $\hat\rho^+, \hat\rho^-$ for $\alpha = 0.95$ ,
and construct two sets of deployment scenarios using \texttt{EXT} and $\texttt{CCG}$.
To assess out-of-sample performance, 
we sample an additional $1\,000$ scenarios (representing realized errors), 
and solve the RT problem \eqref{rt-prob}, 
for $\vc^{\textrm{viol}}=1\,000 \,\textrm{\$/MWh}$.
In addition to \texttt{DSW}, \texttt{EXT}, and $\texttt{CCG}$, 
we solve \eqref{2ro-reform} to optimality via vertex enumeration (\texttt{V-enum}).

\begin{figure}[tb!]
\centerline{\includegraphics[width = \columnwidth]{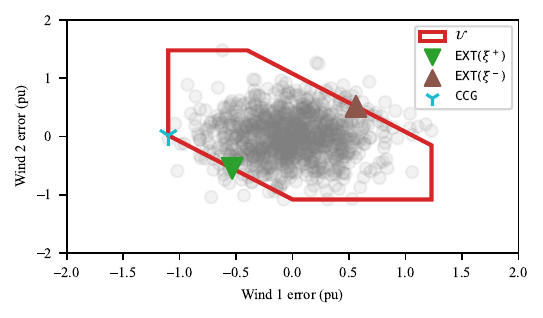}}
\caption{Deployment scenarios found for the 5-bus system ($\alpha=0.95$). 
Grey points indicate the sampled scenarios, \black{used to construct the uncertainty set.}}
\label{5bus_scen_example}
\end{figure}

Fig.~\ref{5bus_scen_example} plots the uncertainty set $\mathcal{U}$, which is \black{enclosed by the} red \black{curve}, and the respective deployment scenarios.
$\texttt{CCG}_{}$ terminates after adding a single scenario, $\vxi^{\textrm{wc}}_1$, 
located at the bottom left vertex, to the DA problem 
whereas \texttt{EXT} adds two deployment scenarios,
and \texttt{V-enum} adds all  six extreme points shown in Fig.~\ref{5bus_scen_example}.
Table~\ref{table_ill_example} presents the average results \black{for scenarios with} realized error \black{that belongs to }$\mathcal{U}$.
Both \texttt{EXT} and $\texttt{CCG}_{}$ lead to higher DA costs and lower RT penalties compared to \texttt{DSW}.
For \texttt{EXT}, $3.3\%$ of out-of-sample observations within $\mathcal{U}$ \black{violate RT constraints (i.e., have at least one RT slack variable with \black{non-zero} value)}. This is 
 an order of magnitude lower than the $54.4\%$ obtained under \texttt{DSW}.
$\texttt{CCG}_{}$ performs best in terms of reliability, with $0\%$ violations. \black{In this case, $\texttt{CCG}_{}$ found the same solution as} $\texttt{V-enum}$, while only using a single scenario.

Examining the DA decisions \black{obtained} by \texttt{DSW}, 
we observe that the generator at node $5$ is the sole reserve supplier, line $1\textrm{-}5$ is congested, and lines $1\textrm{-}4$, $4\textrm{-}5$ are close to becoming congested.
When upward reserves are needed due to negative errors in the wind plant $1$ (over-forecast of wind production), the orange lines shown in Fig.~\ref{5bus_schematic} become congested and the reserves procured at node $5$ cannot be delivered.
The \texttt{EXT} deployment scenario $\vxi^+$ accounts for the impact of constraints on reserve deliverability, but does not consider the worst-case of forecast errors in wind plant $1$ \black{as the grid is ignored during scenario construction}. 
In contrast, the $\texttt{CCG}_{}$ deployment scenario $\vxi^{\textrm{wc}}_1$ fully mitigates reserve undeliverability within the uncertainty set $\mathcal{U}$. 
In the DA solution derived by $\texttt{CCG}_{}$, 
the more expensive generator at node $3$ provides some reserves to compensate for the potential forecast errors in wind plant $1$, 
while the energy schedule remains unchanged.
{While, in this example, the DA schedule changed only in terms of reserve schedules, 
we will see in the next section that inclusion of deployment scenarios in DA can lead to changes in the energy schedule and the available transmission capacity for recourse actions.}

\begin{table}[tb!]
\caption{Results for the 5-bus system ($\alpha=0.95$).}
\centering
\LARGE
\resizebox{\columnwidth}{!}{
\begin{tabular}{@{}llll@{}}
\toprule
    & DA cost ($10^3\textrm{\$/h}$) & Av. RT cost ($10^3\textrm{\$/h}$)  & Viol. Prob. (\%) \\ 
\midrule
\texttt{DSW} & \textbf{131} & 134 & 54.4 \\
\texttt{EXT} & 140 & 1 & 3.3 \\
$\texttt{CCG}_{}$ & 142 & \textbf{0} & \textbf{0} \\ 
\texttt{V-enum} & 142 & \textbf{0} & \textbf{0} \\ 
\bottomrule
\end{tabular}
\label{table_ill_example}}
\end{table}

\subsection{RTS-GMLC 2019 System} \label{large_system_subsection}

\begin{figure}[tb!]
\centering
\includegraphics[width=\columnwidth]{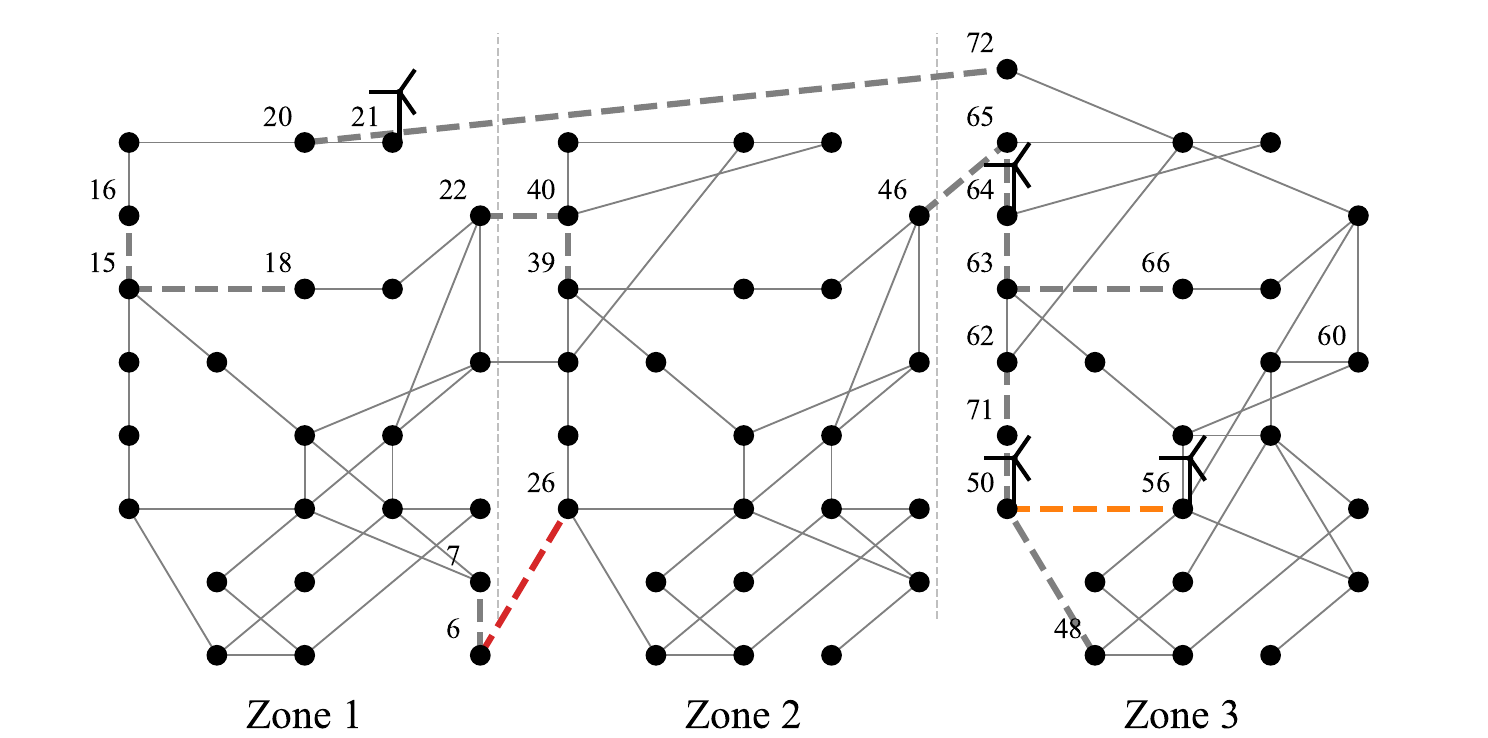}
\caption{Grid topology of RTS-GMLC 2019 System. 
Red color indicates line congestion after \texttt{DSW} is solved and orange color indicates RT line congestion given the \texttt{DSW} decisions, for \black{an illustrative period. Dashed lines are used to initialize ADM.}}
\label{rts_grid_scheme}
\end{figure}

\paragraph{System Information} 
The RTS-GMLC 2019 System \cite{rts_gmlc} has 73 buses, 120 transmission lines, 73 conventional generators, 4 wind power plants, and 56 solar power plants, 
organized in 3 zones (Fig. \ref{rts_grid_scheme}).
We assume perfect demand forecasts and consider imperfect forecasts for the 60 VRE plants.
We use time series data provided by \cite{arpa_e_perform_data} and probabilistic forecasts provided by \cite{carmona2024_joint_granular_model},
comprising $K=500$ scenarios.
The data set covers a full year at hourly granularity.
Exploratory data analysis indicated that the 4 wind plants at nodes $21, 50, 56$, and $64$, and the solar plant at node $60$ account for approximately $75\%$ of the total forecast errors in terms of absolute magnitude.

\paragraph{Implementation Hyper-parameters} 
We compare \texttt{DSW}, \texttt{EXT}, and $\texttt{CCG}_{}$ for levels of $\alpha=\{0.90, 0.95, 0.99\}$, 
\black{ and use a high violation penalty of $\vc^{\textrm{viol}}=1\,000 \,\textrm{\$/MWh}$.}
For the DA scheduling problem \eqref{dsw}, 
we also include an additional variable that allows for VRE curtailment, bounded by the respective point forecast.
For \texttt{CCG}, 
we run the alternating direction method with $L=20$, \black{which is never reached, 
and set the maximum number of scenarios for each period at $M=10$.} 
\black{We initialize Algorithm~\ref{ADM_algo} with the extreme scenarios derived for \texttt{EXT} and with scenarios found using \eqref{adm_init} for $15$ lines that are during at least one period congested in DA or RT, 
according to the \texttt{DSW}-based results --- see Fig. \ref{rts_grid_scheme} for a visualization of the selected lines}.

\begin{table}[tb!]
\caption{Results for a simulation year for the RTS-GMLC system.}
\centering
\LARGE
\resizebox{\columnwidth}{!}{
\begin{tabular}{@{}llccc@{}}
\toprule
 & & Av. DA cost (10$^3\$/\textrm{h})$ & Av. RT cost (10$^3\$/\textrm{h})$ & Viol. Prob. (\%) \\
\midrule
$\alpha=0.90$ & \texttt{DSW} & \textbf{64.57} & 10.69 & 16.96 \\
& \texttt{EXT} & 65.65 & 1.00 & 2.92 \\
& $\texttt{CCG}_{}$ & 78.31 & \textbf{0.02} & \textbf{0.08} \\ 
\midrule
$\alpha=0.95$ & \texttt{DSW} & \textbf{68.24} & 11.20 & 16.92 \\
& \texttt{EXT} & 69.66 & 0.94 & 2.40 \\
& $\texttt{CCG}_{}$ & 82.51 & \textbf{0.01} & \textbf{0.15} \\
\midrule
$\alpha=0.99$ & \texttt{DSW} & \textbf{75.91} & 12.00 & 17.14 \\
& \texttt{EXT} & 78.53 & 0.51 & 2.56 \\
& $\texttt{CCG}_{}$ & 90.26 & \textbf{0.06} & \textbf{0.55} \\
\bottomrule
\end{tabular}
\label{table_RTS}}
\end{table}

\begin{figure}[tb!]
\centering
\includegraphics[width=\columnwidth]{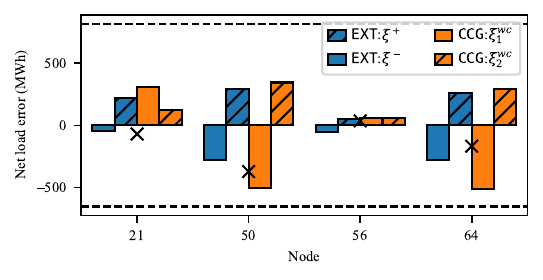}
\caption{Deployment scenarios for $\texttt{EXT}$ and $\texttt{CCG}$ \black{for an illustrative period}.
Dashed lines indicate $\hat\rho^+$, $\hat\rho^-$. 
The $\times$ marker indicates realized errors.}
\label{error_barplot}
\end{figure}

\paragraph{Annual Performance}
Table~\ref{table_RTS} presents annual averages for performance metrics under the different approaches, for $\alpha=\{0.90, 0.95, 0.99\}$. 
The metrics are reported for periods with realized RT error falling within the respective $\mathcal{U}$.
{The percentage of observations falling in $\mathcal{U}$ for each experiment is approximately
$71\%$ for $\alpha = 0.90$; $76\%$ for $\alpha = 0.95$;
and $82\%$ for $\alpha = 0.99$.}
\black{
This undercoverage is attributed to the quality of the probabilistic forecasts, which are an exogenous input in this work.}
From Table~\ref{table_RTS}, as $\alpha$ increases, the DA cost increases and RT cost decreases for all methods, which is expected as the reserve requirements increase with $\alpha$. 
\texttt{DSW}, as expected, performs worst in terms of reliability, 
with an average violation frequency of approximately $17\%$ for all $\alpha$, but performs best in terms of DA cost.
{Considering that we report the frequency of violations when the realized error falls in $\mathcal{U}$, 
we would expect a fully reliable method to have a frequency of violations close to $0$, i.e., assuming \eqref{2ro-prob} is solved to optimality.}
Indeed, \texttt{CCG} achieves a violation frequency $\leq1\%$ in all cases,
whereas \texttt{EXT} is second-best in terms of reliability. 
\black{To understand why these few RT violations within $\mathcal{U}$ persist, 
we examine whether the DA schedule is feasible without using any slack variables (i.e., $\eta^*=0$). For instance, for $\alpha=0.95$, $10$ simulated periods have non-zero RT penalty cost and belong to $\mathcal{U}$. Among these $10$ periods, $8$ have $\eta^*=0$ \footnote{\black{Two periods have $\eta^*>0$, which either indicates that $c^{\textrm{viol}}$ is low or that there is no feasible DA schedule that guarantees zero constraint violations in $\mathcal{U}$. In practice, for such cases, the operator can increase the value of $c^{\textrm{viol}}$.}}. The RT violations during these periods indicate that the ADM failed to find the worst-case error scenario.}

\black{Whereas \texttt{CCG} does not provide any guarantees about reliability outside $\mathcal{U}$, the results are similar to the ones inside $\mathcal{U}$, 
with \texttt{CCG} leading to a lower violation frequency and higher DA cost compared to \texttt{EXT}.
\black{For example, for $\alpha=0.95$, the actual error is outside $\mathcal{U}$ for $24\%$ of the periods. Over these periods, \texttt{EXT} and \texttt{CCG} have an average DA cost of $66.85$ and $80.99 \; 10^{3}\$/\textrm{h}$, respectively; while they record RT violations for $52.11 \%$ and $39.65 \%$ of the periods, respectively.}}

\black{The \texttt{CCG} achieves better performance in terms of reliability by adjusting the choice and number of vertices added as deployment scenarios, considering the dynamic congestion patterns.}
\black{For instance, for $\alpha=0.95$, $\texttt{CCG}_{}$ recovers $3$ or more deployment scenarios in approximately $36\%$ of the time.}
\black{Lastly, it is worth noting that while \texttt{EXT} procures the same amount of upward and downward reserve as \texttt{DSW},  this is not the case for \texttt{CCG} which might choose to procure additional reserve due to congestion anticipated in deployment scenarios with spatial distribution of forecast errors different than \texttt{EXT}. 
For $\alpha=0.95$, in respectively $10.9\%$ and $20.8\%$ of the periods \texttt{CCG} procures more up or down reserve than \texttt{DSW}, 
which implies that the full aggregation benefit cannot be leveraged due to transmission constraints. 
For these periods, on average $3.6$ deployment scenarios are included in the DA schedule.}

\paragraph{Impact of Initialization}
\textcolor{black}{
We examine the sensitivity of the CCG algorithm w.r.t. the starting points by running it with an alternative set of starting points, including only the extreme scenarios.
This results for $\alpha = 0.90/0.95/0.99$ in a violation probability in $\mathcal{U}$ of $1.34/1.36/1.48\%$, which is better than \texttt{EXT} and worse than the results in Table~\ref{table_RTS}, highlighting the benefit of selecting starting points using the initialization strategy in \eqref{adm_init}.}

\paragraph{Illustrative Period}
\black{We illustrate how \texttt{CCG} constructs better deployment scenarios by examining a particular period (with $\alpha = 0.90$) for which the actual error realization falls inside $\mathcal{U}$ and \texttt{EXT} leads to deliverability issues (i.e., non-zero RT penalty), while \texttt{CCG} does not.}
For the selected period, the only VRE generating electricity is wind;
Fig.~\ref{error_barplot} plots the deployment scenarios of $\texttt{EXT}$ and $\texttt{CCG}$ alongside the system-wide requirements $\hat\rho^+, \hat\rho^-$ and the realized error, 
showing only nodes $21, 50, 56$ and $64$, corresponding to the four wind plants.
\black{For both \texttt{EXT} and \texttt{CCG}, the aggregate error in each deployment scenario is equal to either $\hat\rho^+$ or $\hat\rho^-$, but its distribution among nodes differs.}
For $\texttt{CCG}$, the first deployment scenario, $\vxi^{\textrm{wc}}_1$, 
has an aggregate error equal to $\hat\rho^-$;
compared to $\vxi^-$, $\vxi^{\textrm{wc}}_1$ leads to a much higher allocation in nodes $50, 64$ \black{(RT wind production higher than forecasted)}
whereas nodes $21,56$ have a smaller error with opposite sign; 
thus, $\vxi^{\textrm{wc}}_1$ represents a less correlated setting.
Concerning the second deployment scenario, $\vxi^{\textrm{wc}}_2$, we observe that all nodal errors are positive \black{(RT wind production lower than forecasted)}
with the aggregate error equal to $\hat\rho^+$;
and the allocation is not proportional to the width of the prediction interval at each node, with relatively to $\texttt{EXT}$ a higher allocation at nodes $50, 64$. 
The line connecting nodes $6\textrm{-}26$ is congested in the \texttt{DSW} solution of the DA problem.

\begin{figure}[tb!]
\centerline{\includegraphics[width = \columnwidth]{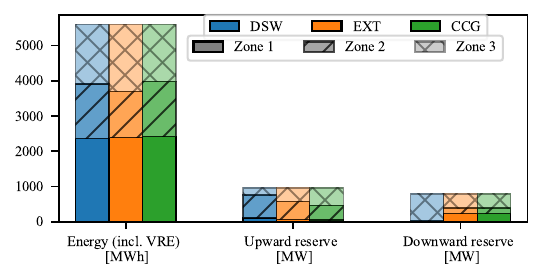}}
\caption{DA schedule for an illustrative period, aggregated per zone.}
\label{stacked_barplot_da_schedule}
\end{figure}

In RT, a large negative error \black{(RT wind production higher than forecasted)}
occurs in node $50$ (see $\times$ marker in Fig.~\ref{error_barplot});
in turn, under \texttt{DSW}, \black{the re-dispatch} causes congestion in line $50\textrm{-}56$ and leads to activating the respective slack variable (approximately $96 \textrm{MW}$);
the same slack is activated in RT for \texttt{EXT} but at a lower level than in  \texttt{DSW} ($47 \textrm{MW}$), 
meaning that \texttt{EXT} only partially resolves the deliverability issue. 
In contrast, \texttt{CCG} does not activate any slacks in RT.
To understand this effect, we examine the aggregate energy and reserve schedule illustrated in Fig.~\ref{stacked_barplot_da_schedule}.
Firstly, \texttt{CCG} schedules more energy in zone 2 by curtailing \black{a portion of wind production} in node $50$, 
thus creating available transmission capacity for recourse actions of $187 \textrm{MW}$ in line $50\textrm{-}56$, which is higher than the {$97 \textrm{MW}$ created by \texttt{EXT}}.
Secondly, \texttt{CCG} shifts upward reserve schedule from zone 1 and 2 to zone 3, 
whereas \black{downward reserves shift from zone 3 to zone 1 and 2}, 
which is due to the fact that the deployment scenarios consider higher forecast errors in nodes $50, 64$.
\black{This can also be understood from the sign of the respective entries within the PTDF matrix, 
which are negative for nodes $50, 64$ towards line $50-56$ and thus the second \texttt{CCG}-scenario creates bigger available transmission capacity for recourse actions than \texttt{EXT}.}
The rest of the time periods offer similar insights.
Namely, the $\texttt{CCG}$ deployment scenarios allocate the uncertainty non-proportionally across nodes.
As a result, $\texttt{CCG}$ \black{yields adapted energy-reserve schedules and available transmission capacity for recourse actions,} fit for balancing different spatial distributions of forecast errors.

\section{Conclusion}\label{section_conclusions}

Transmission congestion often renders reserves undeliverable during real-time operations, threatening system reliability.
Emerging industry practices add a set of deployment scenarios as contingencies to the day-ahead \black{scheduling problem} to enhance reserve deliverability during real-time operations. 
However, industry practices for constructing deployment scenarios typically ignore the interdependency between deployment scenarios, grid characteristics, and day-ahead schedules.
In this work, we leverage adaptive robust optimization to formulate a two-stage problem that jointly considers the day-ahead schedule and worst-case scenarios for reserve deliverability, and solve it via a column-and-constraint algorithm to generate reserve
deployment scenarios.
\black{We conduct simulations for the RTS-GMLC 2019 system using two sets of scenarios: one set is constructed via the CCG algorithm and another based on a method that resembles current industry practice. } 
Overall, the deployment scenarios constructed via the CCG algorithm significantly lower the \black{number of periods} with \black{undeliverable reserves} and are more appropriate to meet the system operator's reliability targets.
\black{Contrary to prevailing industry heuristics}, these deployment scenarios often include forecast errors in opposite directions, \black{which have the effect of aggravating the congestion of transmission lines}.
In addition, simulation results show that the number of deployment scenarios and the spatial distributions of errors dynamically adapt to the day-ahead conditions.
\textcolor{black}{Future work could quantify the benefits of the proposed approach for scheduling processes that include commitment decisions and inter-temporal constraints.}
\black{Further work could expand the CCG algorithm to include provisions for cases where no feasible day-ahead schedule is found to be robust within the uncertainty set. Alternative algorithms for solving the adversarial problem could also be tested.}
Lastly, it is interesting to examine generating deployment scenarios directly from data, thus bypassing the need for iterative solution techniques.

\bibliographystyle{IEEEtran}
\bibliography{IEEEabrv,references}

\end{document}